\documentclass[a4paper,10pt]{scrartcl}
\usepackage[english]{babel}
\usepackage[utf8]{inputenc}
\usepackage[T1]{fontenc}
\usepackage{lmodern} 
\usepackage{enumerate}
\usepackage[a4paper]{geometry}
\usepackage{amsfonts,amsmath,amssymb,amsthm,amsopn}
\usepackage{mathtools}
\usepackage{stmaryrd}
\usepackage{nicefrac}
\usepackage{exscale}

\usepackage[numbers,square]{natbib}
\usepackage{url} 
\usepackage[unicode]{hyperref}

\usepackage{graphicx} 
\usepackage{scrpage2}

\allowdisplaybreaks 

\newtheoremstyle{theo}
	{3pt} 
	{3pt} 
	{\itshape} 
	{} 
		{\bfseries} 
	{\\} 
	{ } 
	{\thmname{#1}\thmnumber{ #2.}\thmnote{ - #3}} 
\theoremstyle{theo}

\newtheorem{defi}{Definition}[section]
\newtheorem{lem}[defi]{Lemma}
\newtheorem{theo}[defi]{Theorem}

\newtheorem{rem}[defi]{Remark}
\newtheorem{prop}[defi]{Proposition}

\newenvironment{bew}{\begin{proof}[\bfseries Proof:]}{\end{proof}}

\DeclareMathOperator{\bomega}{\overline{\Omega}}
\DeclareMathOperator{\romega}{\partial\Omega}
 
\DeclareMathOperator{\Lo}{L}
\DeclareMathOperator{\W}{W}

\DeclareMathOperator{\intd}{d\!}
\DeclareMathOperator{\divdot}{\!\cdot}
\DeclareMathOperator{\ep}{(\varepsilon)}
\DeclareMathOperator{\chie}{\chi^{(\varepsilon)}}

\newcommand{\intos}{\int\limits_0^{s_0}}

\newcommand{\sintomega}{\smallint_{\text{\tiny{$\Omega$}}}} 
\newcommand{\sintRn}{\smallint_{\text{\tiny{$\mathbb{R}^n$}}}} 

\newcommand{\Wu}{\underline{W}} 
\newcommand{\R}{\mathbb{R}}

\newcounter{gleichung}
\newcommand{\owncount}{\refstepcounter{gleichung}}

\author{Tobias Black\thanks{Institut f\"ur Mathematik, Universit\"at Paderborn, Warburger Str. 100, 33098 Paderborn, Germany; email: \mbox{tblack@math.upb.de}}}
\title{Blow-up of weak solutions to a chemotaxis system under influence of an external chemoattractant}

\begin{document}
\maketitle

\begin{abstract}
\noindent
{\textbf{Abstract: }
We study nonnnegative radially symmetric solutions of the parabolic-elliptic Keller-Segel whole space system
\begin{align*}
\left\{\begin{array}{c@{\,}l@{\quad}l@{\,}c}
u_{t}&=\Delta u-\nabla\divdot(u\nabla v),\ &x\in\mathbb{R}^n,& t>0,\\
0 &=\Delta v+u+f(x),\ &x\in\mathbb{R}^n,& t>0,\\
u(x,0)&=u_{0}(x),\ &x\in\mathbb{R}^n,&
        \end{array}\right.
\end{align*}
with prototypical external signal production
\begin{align*}
f(x):=\begin{cases}
f_0\vert x\vert^{-\alpha},&\text{ if }\vert x\vert \leq R-\rho,\\
0,&\text{ if } \vert x\vert\geq R+\rho,\\
\end{cases}
\end{align*}
for $R\in(0,1)$ and $\rho\in\left(0,\frac{R}{2}\right)$, which is still integrable but not of class $\text{L}^{\frac{n}{2}+\delta_0}(\mathbb{R}^n)$ for some $\delta_0\in[0,1)$. 
For corresponding parabolic-parabolic Neumann-type boundary-value problems in bounded domains $\Omega$, where $f\in\text{L}^{\frac{n}{2}+\delta_0}(\Omega)\cap C^{\alpha}(\Omega)$ for some $\delta_0\in(0,1)$ and $\alpha\in(0,1)$, it is known that the system does not emit blow-up solutions if the quantities $\|u_0\|_{\text{L}^{\frac{n}{2}+\delta_0}(\Omega)}, \|f\|_{\text{L}^{\frac{n}{2}+\delta_0}(\Omega)}$ and $\|v_0\|_{\text{L}^{\theta}(\Omega)}$, for some $\theta>n$, are all bounded by some $\varepsilon>0$ small enough.

We will show that whenever $f_0>\frac{2n}{\alpha}(n-2)(n-\alpha)$ and $u_0\equiv c_0>0$ in $\overline{B_1(0)}$, a measure-valued global-in-time weak solution to the system above can be constructed which blows up immediately. Since these conditions are independent of $R\in(0,1)$ and $c_0>0$, we will thus prove the criticality of $\delta_0=0$ for the existence of global bounded solutions under a smallness conditions as described above.
}
\end{abstract}
\pagebreak

\section{Introduction}\label{ch:c1}
Chemotaxis is a biological mechanism whereby the movement of cells is influenced by a chemical substance. This mechanism appears in multiple biological processes, e.g. aggregation of bacteria or the inflammatory response of leukocytes. One of the first PDE systems modelling these processes dates back to the pioneering works \cite{KS70} and \cite{KS71} by Keller and Segel. Variants of the original Keller-Segel model have also been incorporated in more complex biological processes ranging from pattern formation (\cite{AMM10}) to agiogenesis in early stages of cancer (\cite{SRLC09}).
For a broader spectrum of applications and an overview of known results we refer to the survey articles \cite{BBWT15}, \cite{HP09} and \cite{Ho03}.

The Keller-Segel system which is the basis this article has the form
\begin{align}\label{KS}\tag{$KS$}
\left\{\begin{array}{c@{\,}l@{\quad}l@{\,}c}
u_{t}&=\Delta u-\nabla\divdot(u\nabla v),\ &x\in\Omega,& t>0,\\
v_{t}&=\Delta v-v+u,\ &x\in\Omega,& t>0,\\
\frac{\partial u}{\partial\nu}&=\frac{\partial v}{\partial\nu}=0,\ &x\in\romega,& t>0,\\
 u(x,0)&=u_{0}(x),\ v(x,0)=v_{0}(x),\ &x\in\Omega,&
        \end{array}\right.
\end{align}
wherein $u(x,t)$ represents the density of the moving cells and $v(x,t)$ denotes the concentration of an attracting chemical substance influencing said movement at place $x$ in the bounded domain $\Omega$ and at time $t$. In the mathematical study of chemotaxis, blow-up solutions, i.e. the existence of some $T\in(0,\infty]$ such that $\limsup_{t\nearrow T}\|u\|_{\Lo^\infty(\Omega)}=\infty$, are of utmost importance. The existence of such solutions is identified with the occurrence of self-organizing patterns within the cell population. As such, the formulation of conditions which allow for blow-up to happen, or conditions negating blow-up completely are widely sought after.

For \eqref{KS} conditions negating blow-up are well known. In particular it was shown for suitable $\Omega\subset\R^n$ and nonnegative initial values $u_0\in C^0(\bomega)$, $v_0\in C^1(\bomega)$, that the corresponding maximally extended classical solution $(u,v)$ of \eqref{KS} fulfills:

\begin{itemize}
\item[If] $n=2:$ If $\sintomega u_0\intd x<4\pi$ (or $8\pi$ in the radial symmetric setting), then $(u,v)$ is global and bounded with regard to the $\Lo^\infty(\Omega)$-norm. (\cite{NSY97}, \cite{gz98})
\item[If] $n\geq3:$ It was proven in \cite{win10jde} that there exists a bound for $u_0$ in $\Lo^q(\Omega)$ and for $\nabla v_0$ in $\Lo^p(\Omega)$, with $q>\frac{n}{2}$ and $p>n$ such that the solution $(u,v)$ is global in time and bounded. This result has further been extended to the critical case $q=\frac{n}{2}$ and $p=n$. (\cite{cao2014global})
\end{itemize}

In our previous work (\cite{Bla15}) we considered an extension of the \eqref{KS} model by introducing an external signal production to \eqref{KS}. Namely, we studied the system
\begin{align}\label{KSf}\tag{$KS_f$}
\left\{\begin{array}{c@{\,}l@{\quad}l@{\,}c}
u_{t}&=\Delta u-\nabla\divdot(u\nabla v),\ &x\in\Omega,& t>0,\\
\tau v_{t}&=\Delta v-v+u+f(x,t),\ &x\in\Omega,& t>0,\\
\frac{\partial u}{\partial\nu}&=\frac{\partial v}{\partial\nu}=0,\ &x\in\romega,& t>0,\\
 u(x,0)&=u_{0}(x),\ v(x,0)=v_{0}(x),\ &x\in\Omega&
        \end{array}\right.
\end{align}
in a bounded and smooth domain $\Omega\subset\R^n$ with $n\geq2,\tau>0,u_0\in C^0(\bomega), v_0\in W^{1,\theta}(\Omega)$ for some $\theta>n$ and $f\in\Lo^\infty\left([0,\infty);\Lo^{\frac{n}{2}+\delta_0}(\Omega)\right)\cap C^\alpha\left(\Omega\times(0,\infty)\right)$ with $\delta_0\in(0,1)$ and $\alpha>0$. Obviously for $f\equiv0$, this coincides  with $\eqref{KS}$. 
Applying the usual fixed point arguments, as illustrated in \cite{win10ctax} for the system \eqref{KS}, we were able to verify the existence of classical solutions to \eqref{KSf}. Furthermore, we were able to expand most of the known boundedness results for \eqref{KS} to the setting $f\not\equiv 0$. Firstly, we obtained a critical mass result for $n=2$ and constant signal production $f\in\Lo^{\frac{n}{2}+\delta_0}(\Omega)\cap C^\alpha\left(\Omega\right)$ (\cite[Theorem 1.2]{Bla15}) similar to the one by \cite{NSY97} and \cite{gz98}. Secondly, and important for the context of our current work, we were able to prove a result similar to the one above by \cite{win10jde}. We cite the theorem here in a short version without including the statements regarding asymptotic behavior of the solution. For the full version see \cite[Theorem 1.3]{Bla15}.

\begin{theo}
Let $0<\delta_0<1$, $0<\alpha$, $n<\theta<\frac{n^2+2n\delta_0}{n-2\delta_0}$ and $1<r$. Then there exist constants $\varepsilon_0>0$ and $C>0$ with the following property:
If $u_0\in C^0(\bomega)$, $v_0\in\W^{1,\theta}(\Omega)$ and $f\in\Lo^\infty([0,\infty);\Lo^{\frac{n}{2}+\delta_0}(\Omega))\cap C^\alpha(\Omega\times(0,\infty))$ are nonnegative with
\begin{align*}
\|u_0\|_{\Lo^{\frac{n}{2}+\delta_0}(\Omega)}\leq\varepsilon,\ \|\nabla v_0\|_{\Lo^\theta(\Omega)}\leq\varepsilon\text{ and }\|f\|_{\Lo^\infty([0,\infty);\,\Lo^{\frac{n}{2}+\delta_0}(\Omega))}\leq\varepsilon
\end{align*}
for some $\varepsilon<\varepsilon_0$, then there exists a global classical solution $(u,v)$ of \eqref{KSf} with $\|u\|_{\Lo^\infty(\Omega)}$ and $\|v\|_{\W^{1,\theta}(\Omega)}$ remaining bounded for all times. 
\end{theo}

Unfortunately, the methods applied to prove the result above do not yield any information whether $\delta_0=0$ is the critical boundary for the existence of such small-data solutions. This seems to be strongly suggested by the results from \cite{TW11}, where a simplified parabolic-elliptic version of \eqref{KSf} in the radially symmetric setting on the whole space $\mathbb{R}^n$ for $n=2$ was considered, that is 
\begin{align}\label{KS0f}\tag{$KS^{\,0}_f$}
\left\{\begin{array}{c@{\,}l@{\quad}l@{\,}c}
u_{t}&=\Delta u-\nabla\divdot(u\nabla v),\ &x\in\mathbb{R}^n,& t>0,\\
0 &=\Delta v+u+f(x),\ &x\in\mathbb{R}^n,& t>0,\\
u(x,0)&=u_{0}(x),\ &x\in\mathbb{R}^n,&
        \end{array}\right.
\end{align}
with a Dirac-distributed signal production $f(x)=f_0\delta(x)$. It was shown in a radially symmetric setting that for any choice of $f_0>0$ certain generalized solutions, so called radial weak solutions, blow up immediately and depending on the size of the initial mass $\mu:=\sintRn u_0(x)\intd x<\infty$, compared to the critical mass $8\pi-2f_0$, form a Dirac singularity.

It is the purpose of the present work to examine whether $\delta_0=0$ is also a critical boundary for the existence of such small-data solutions in higher dimensions. To this end we will study the behavior of solutions in dimensions $n\geq3$ for constant-in-time functions $f$, that are not of class $\Lo^{\frac{n}{2}+\delta_0}(\R^n)$ for some $\delta_0\in[0,1)$ but still integrable. Following the approach of \cite{TW11}, we will consider \eqref{KS0f} in the radially symmetric setting on the whole space $\mathbb{R}^n$ for $n\geq3$, with given radially symmetric and nonnegative $u_0$. Furthermore we assume, that $u_0\not\equiv 0$ has finite mass $\mu$ and that $f$ is nonnegative and radially symmetric as well.

Using a transformation introduced in \cite{JL92} and employed in \cite{TW11}, we will prove that generalized global-in-time measure-valued solutions of \eqref{KS0f} blow up immediately for prototypical signal production functions $f$ satisfying
\begin{align*}
f(x):=\begin{cases}
f_0\vert x\vert^{-\alpha},&\text{ if }\vert x\vert \leq R-\rho,\\
0,&\text{ if } \vert x\vert\geq R+\rho\\
\end{cases}
\end{align*}
and smooth in between with some $1>R>0,\alpha>2$, $\rho\in(0,\frac{R}{2})$ and $f_0>\frac{2n}{\alpha}(n-2)(n-\alpha)$. 

Our main theorem reads as follows:
\begin{theo}[Immediate blow-up of radial weak solutions]\label{theoblowup}
For $1>R>0,\rho\in(0,\frac{R}{2}), n>\alpha>2$ and $f_0>0$ satisfying $f_0>\frac{2n}{\alpha}(n-2)(n-\alpha)$ let $f(r)$ be defined as in \eqref{fdef}. Furthermore, assume the initial data satisfy $u_0\equiv c_0$ in $\overline{B_1(0)}$  for some positive constant $c_0>0$. Then for all $t_0\geq0$ there exists a globally defined radial weak solution $u$ of \eqref{KS0f}, in the sense of Definition \ref{weakrad} below, such that this solution satisfies
\begin{align*}
\|u\|_{\Lo^\infty(\mathbb{R}^n\times(t_0,t_0+\eta))}=\infty\quad\text{ for all }\eta>0.
\end{align*}
\end{theo}
The theorem above states a sufficient condition for the occurrence of immediate blow-up in \eqref{KS0f}, with external production of the form described in \eqref{fdef}. The only restriction on $u_0\in\Lo^1(\R^n)\cap\Lo^\infty(\R^n)$ involves $c_0>0$, which can be arbitrary small. Obviously, one can thereby find smooth initial values with arbitrary small $\Lo^p$ norms for which the theorem is still applicable and since the assumption imposed on $f_0$ is independent of the parameter $R$, even signal production functions with small norm may lead to blow-up -- in the case where $f$ is not of class $\Lo^{\frac{n}{2}+\delta_0}(\mathbb{R}^n)$ (for some $0\leq\delta_0<1$). Thus, the case $\delta_0=0$ is indeed critical. We have to leave open the question if blow-up may also occur for scaling factors $f_0$ smaller than $\frac{2n}{\alpha}(n-2)(n-\alpha)$, since the methods used here give no evidence on the behavior of solutions when the scaling factor of the prototypical signal production is small.

In the subsequent sections, if not stated otherwise, $n$ will always denote the space dimension and $\mu:=\sintRn u_0(x)\intd x<\infty$ the initial mass.

%
\setcounter{gleichung}{0} 
\section{Blow-up for less regular signal production}\label{cmeasval}
We would like to consider prototypical signal production functions of the form
\begin{align}\label{fhdef}
\owncount
\tilde{f}(x):=\begin{cases}
f_0\vert x\vert^{-\alpha},&\text{ if }\vert x\vert < R,\\
0,&\text{ if } \vert x\vert\geq R
\end{cases}
\end{align}
with some $f_0>0$ and $0<R<1$, which for $\alpha\in(2,n)$ are still integrable but not of class $\Lo^{\frac{n}{2}+\delta_0}(\mathbb{R}^n)$ for any $\delta_0\in[0,1)$. However, in order to take advantage of well-known regularity results for partial differential equations we will instead work with a smoother version of \eqref{fhdef}. More precisely, for $\rho\in(0,\frac{R}{2})$ we consider radially symmetric functions $f$ satisfying
\begin{align}\label{fdef}
\owncount
f(r):=\begin{cases}
f_0 r^{-\alpha},&\text{ if } r \leq R-\rho,\\
0,&\text{ if } r\geq R+\rho,
\end{cases}
\end{align}
such that $r\mapsto f(r)$ is smooth and monotonically decreasing on $(0,\infty)$. Consequently the function $F(s):=\smallint\limits_0^{s^{\nicefrac{1}{n}}}f(r)r^{n-1}\intd r$ is monotonically increasing and smooth on $(0,\infty)$ with
\begin{align}\label{Fdef}
\owncount
F(s)=\smallint\limits_0^{s^{\nicefrac{1}{n}}}f(r)r^{n-1}\intd r\begin{cases}
=\frac{f_0}{n-\alpha}s^{\frac{n-\alpha}{n}}&,\text{ if }0\leq s\leq R-\rho,\\
\leq\frac{f_0}{n-\alpha}(R+\rho)^{\frac{n-\alpha}{n}}&,\text{ if }s\geq R+\rho.
\end{cases}
\end{align}
In addition the first derivative is monotonically decreasing and satisfies
\begin{align}\label{Fsdef}
\owncount
F_s(s)=\frac{1}{n}f(s^{\nicefrac{1}{n}})=\begin{cases}
\frac{f_0}{n}s^{-\frac{\alpha}{n}}&,\text{ if }0<s\leq R-\rho,\\
0&,\text{ if }s\geq R+\rho.
\end{cases}
\end{align}
Both of these functions will play an important role in the transformation of \eqref{KS0f} introduced in the next section, which is an adjustment to higher space dimensions of the transformation applied in \cite{TW11}. 
\subsection{Radial weak solutions}\label{cmeasvals1}
Following the approach employed in \cite{TW11} and \cite{JL92}, we will first make use of spherical coordinates to transform \eqref{KS0f} into the related degenerate parabolic initial-boundary value problem \eqref{WDGL}. The solutions of these PDE problems are connected by the notion of radial weak solutions stated in Definition \ref{weakrad}. Our objective is then to prove the immediate blow-up of $W_s$ and in turn, by the identity \eqref{backtrafo}, the blow-up of the radial weak solution $u$.

The transformation in question is defined by
\begin{align}\label{wtrafo1}
\owncount
W(s,t):=\frac{n}{\vert S_{n-1}\vert}\int\limits_{B\left(0,\,s^{\nicefrac{1}{n}}\right)}u(x,t)\intd x,\ s\geq0,t\geq0,
\end{align}
with $\vert S_{n-1}\vert$ representing the surface area of the unit sphere in $n$ dimensions and $B(0,r)$ denoting the ball around the origin with radius $r$. For radially symmetric $u=u(r,t)$, by using spherical coordinates, this can also be expressed as
\begin{align*}
W(s,t)=\frac{n}{\vert S_{n-1}\vert}\int\limits_{S_{n-1}}\int\limits_0^{s^{\nicefrac{1}{n}}}u(r,t)r^{n-1}\intd r\intd S=n\int\limits_0^{s^{\nicefrac{1}{n}}}u(r,t)r^{n-1}\intd r,\ s\geq0,t\geq0.
\end{align*}
Formal calculation, without regarding the regularity of $u$ for now, shows
\begin{align}\label{wtrafo2}
\owncount
W_s(s,t)=n u(s^{\frac{1}{n}},t)(s^{\frac{1}{n}})^{n-1}\frac{1}{n}s^{\frac{1-n}{n}}=u(s^{\frac{1}{n}},t)\text{ for } s>0,t\geq0
\end{align}
and thus
\begin{align*}
W_{ss}(s,t)=\frac{1}{n}s^{\frac{1-n}{n}}u_r(s^{\frac{1}{n}},t)\text{ for } s>0,t\geq0.
\end{align*}
Considering these expressions and the first equation of \eqref{KS0f} we thereby see that $W(s,t)$ formally fulfills
\begin{align*}
W_t(s,t)&=\frac{n}{\vert S_{n-1}\vert}\!\int\limits_{B\left(0,\,s^{\nicefrac{1}{n}}\right)}\!\!\!\!\Delta u-\nabla u\nabla v-u\Delta v\intd x
=n\int\limits_0^{s^{\nicefrac{1}{n}}}(r^{n-1}u_r)_r-(u r^{n-1}v_r)_r\intd r\\
&=n\left(s^{\frac{n-1}{n}}u_r(s^{\nicefrac{1}{n}},t)-\left[u v_r r^{n-1}\right]^{s^{\nicefrac{1}{n}}}_0\right).
\end{align*}
The second equation of \eqref{KS0f} implies $r^{n-1}v_r=-\int\limits_0^r\sigma^{n-1}u(\sigma,t)\intd\sigma-\int\limits_0^r\sigma^{n-1}f(\sigma)\intd\sigma$, and thus
\begin{align*}W_t(s,t)&=n^2 s^{\frac{2n-2}{n}}W_{ss}(s,t)+W(s,t)W_s(s,t)+nF(s)W_s(s,t)
\end{align*}
with $F(s)$ as in \eqref{Fdef}. Let us briefly recall further statements from \cite{TW11} regarding the transformed problem.
Letting $W_0:=W(s,0)=n\smallint\limits_0^{s^{\nicefrac{1}{n}}}u_0(r)r^{n-1}\intd r\text{ for } s\geq0,$ we observe that if $u_0$ is nonnegative and bounded fulfilling $\mu:=\sintRn u_0(x)\intd x<\infty$, then $W_0$ satisfies
\begin{align}\label{w0prop}
\owncount
W_0&\in\W^{1,\infty}((0,\infty)),\nonumber\\
W_{0s}&\geq0\text{ in }(0,\infty),\\
W_0(s)&\rightarrow\frac{n\mu}{\vert S_{n-1}\vert}\text{ as }s\rightarrow\infty\nonumber.
\end{align}
If $u$ satisfies the mass-conservation property $\sintRn u(x,t)\intd x=\mu$ for all $t\geq0$, then for each $t\geq0$ there holds $W(s,t)\rightarrow\frac{n\mu}{\vert S_{n-1}\vert}$ as $s\rightarrow\infty$.
Thereby, this formally leads to the following degenerate parabolic initial-boundary value problem:
\begin{align}\label{WDGL}
\owncount
\left\{\begin{array}{c@{\,}l@{\quad}l@{\,}c}
W_{t}&=n^2s^{\frac{2n-2}{n}}W_{ss}+W_sW+nF W_s\ &s>0,& t>0,\\
W(0,t) &=0,\ \displaystyle \lim_{s\rightarrow\infty}W(s,t)=\frac{n\mu}{\vert S_{n-1}\vert} ,& t>0,\\
W(s,0)&=W_{0}(s)\ &s>0.&
        \end{array}\right.
\end{align}
By the definition of the transformation in \eqref{wtrafo1} and the nonnegativity of $u$, it is obvious that for each $t$ the function $W(\cdot,t)$ must be nondecreasing. In particular, $W$ is bounded and since the PDE in \eqref{WDGL} is uniformly parabolic with smooth coefficients in each cylinder $(s_0,\infty)\times(0,\infty)$ with $s_0>0$, standard theory implies the smoothness of $W_s$ in $(0,\infty)\times(0,\infty)$. Thus, for $s>0$ we expect the identity $W_s(s,t)=u(s^{\nicefrac{1}{n}},t)$ suggested by \eqref{wtrafo2} to hold. Consequently, discontinuities for $W$ can only occur at $s=0$, and given a solution $W$ of \eqref{WDGL} we can reconstruct $u$ by taking into account the jump size $W(0+,t):=\lim_{s\searrow0}W(s,t)$ at the origin, in terms of the measure-valued identity
\begin{align}\label{backtrafo}
\owncount
u(x,t):=W_s(\vert x\vert^{n},t)+\frac{\vert S_{n-1}\vert}{n} W(0+,t)\delta(x)
\end{align}
for $t>0$, where $\delta(x)$ denotes the Dirac delta function in $n$ dimensions. In \cite{TW11} it was therefore suggested to act in the framework of radially symmetric Radon measures $\mathcal{M}_{rad}(\mathbb{R}^n)$, that is the space of all functionals $\psi$, radially symmetric about $x=0$, defined on the space $C_0^0(\mathbb{R}^n)$ of compactly supported continuous functions over $\R^n$. To be more precise
\begin{align*}
\mathcal{M}_{rad}(\mathbb{R}^n):=\left\{\psi:C_0^0(\R^n)\rightarrow\R\,\right\vert \,\psi\!\left(\zeta\circ\upsilon\right)\!=\psi\left(\zeta\right)&\text{ for all test functions }\zeta\in C^{\infty}_0(\R^n)\\
&\left.\!\text{ and all rotations }\upsilon\in SO(n)\right\},
\end{align*}
where $SO(n)$ denotes the special orthogonal group in $n$ dimensions. This way, we translate the notion of radial weak solutions given in \cite[Definition 1.1]{TW11} to our equation in the following way:
\begin{defi}\label{weakrad}
Assume $u_0\in\Lo^{\infty}(\mathbb{R}^n)$ is nonnegative and $\mu:=\sintRn u_0(x)\intd x$ is finite. Then we call
\begin{align*}
u\in C^0([0,\infty);\mathcal{M}_{rad}(\mathbb{R}^n))
\end{align*}
a radial weak solution of \eqref{KS0f} in $\mathbb{R}^n\times(0,\infty)$ if the function $W:[0,\infty)\times[0,\infty)$ defined by \eqref{wtrafo1} satisfies
\begin{align*}
W(s,t)\rightarrow\frac{n\mu}{\vert S_{n-1}\vert}\ \text{ as }s\rightarrow\infty\ \text{ for all }t>0
\end{align*}
and
\begin{align*}
-\int\limits_0^\infty\int\limits_0^\infty\zeta_t W - \int\limits_0^\infty \zeta(\cdot,0)W_0
=\ n^2\int\limits_0^\infty\int\limits_0^\infty (s^{\frac{2n-2}{n}}\zeta)_{ss}W-\frac{1}{2}\int\limits_0^\infty\int\limits_0^\infty\zeta_s W^2-n\int\limits_0^\infty\int\limits_0^\infty(F\zeta)_s W
\end{align*}
for all $\zeta\in C_0^{\infty}([0,\infty)\times[0,\infty))$, where $W_0(s):=\frac{n}{\vert S_{n-1}\vert}\smallint_\text{\tiny{$B(0,\!s^{\nicefrac{1}{n}})$}}u_0(x)\intd x$ for $s\geq 0$.
\end{defi}
The construction of such solutions is based on an approximation procedure wherein we utilize cut-off functions to counteract the degeneracy of the second and third term on the right side of \eqref{WDGL}. Fortunately cut-off functions of the type used in \cite{TW11} can also be applied in the current setting. Following the approach of \cite{TW11}, we first fix a cut-off function $\chi\in C^{\infty}([0,\infty))$ such that $\chi\equiv0$ on $[0,\frac{1}{2}], \chi\equiv1$ on $[1,\infty)$ and $\chi'\geq0$ on $[0,\infty)$. Then for $\varepsilon\in(0,1)$ we define $\chie(s):=\chi(\frac{s}{\varepsilon})$. This way $\chie$ has the properties
\begin{align}\label{cut-off1}
\owncount
\chie\equiv0\text{ on }[0,\frac{\varepsilon}{2}],\quad \chie\equiv1\text{ on }[\varepsilon,\infty)\quad\text{ and }\quad\chi^{\ep}_s\geq0\text{ on }[0,\infty).
\end{align}
Additionally $\chie$ satisfies the inequalities
\begin{align}
\owncount\label{cut-off2}
\vert\chi^{\ep}_s\vert\leq\frac{c_\chi}{\varepsilon}\quad\text{ and }\quad\vert\chi^{\ep}_{ss}\vert\leq\frac{c_\chi}{\varepsilon^2}
\end{align}
with $c_\chi:=\|\chi'\|_{\Lo^{\infty}((0,\infty))}+\|\chi''\|_{\Lo^\infty((0,\infty))}$. Moreover, $\chie(s)\nearrow1$ as $\varepsilon\rightarrow0$ holds for all $s>0$.

The cut-off function at hand, we now introduce the approximate problem for \eqref{WDGL}:
\begin{align}\label{WDGLeps}
\owncount
\left\{\begin{array}{c@{\,}l@{\quad}l@{\,}c}
W^{\ep}_{t}&=n^2s^{\frac{2n-2}{n}}W^{\ep}_{ss}+\chie W^{\ep}_sW^{\ep}+nF\chie W^{\ep}_s,\ &s>0,& t>0,\\
W^{\ep}(0,t) &=0,\ \displaystyle \lim_{s\rightarrow\infty}W^{\ep}(s,t)=\frac{n\mu}{\vert S_{n-1}\vert} ,& t>0,\\
W^{\ep}(s,0)&=W_{0}(s),\ &s>0.&
        \end{array}\right.
\end{align}
Although we are more interested in the behavior of solutions, we cannot completely skip examining solvability and other important properties. Let us therefore briefly state some results whose proofs we omit since these results can be proven by using the same arguments as shown in \cite[Lemma 1.2 - Lemma 1.5]{TW11}.
\begin{lem}\label{compWeps}
For $1>R>0,\rho\in(0,\frac{R}{2}), n>\alpha>2$ and $f_0>0$ let $F$ be defined as in \eqref{Fdef}. Assume $\overline{W}$ and $\underline{W}$ belong to $C^0([0,\infty)\times[0,\infty))\cap C^{2,1}((0,\infty)\times(0,\infty))$ and satisfy
\begin{align*}
\overline{W}_t&\geq n^2s^{\frac{2n-2}{n}}\overline{W}_{ss}+\chie\overline{W}\overline{W}_s+n\chie F\overline{W}_s
\intertext{and}
\underline{W}_t&\leq n^2s^{\frac{2n-2}{n}}\underline{W}_{ss}+\chie\underline{W}\underline{W}_s+n\chie F\underline{W}_s
\end{align*}
for all $s>0$ and $t>0$. Additionally, suppose $\overline{W}_0:=\overline{W}(\cdot,0)$ and $\underline{W}_0:=\underline{W}(\cdot,0)$ satisfy \eqref{w0prop} with positive numbers $\overline{\mu}$ and $\underline{\mu}$ respectively and that $\overline{W}_0\geq\underline{W}_0$ holds on $(0,\infty)$. Moreover, assume $\overline{W}(0,t)\geq\underline{W}(0,t)$ for all $t>0$ and $\displaystyle\lim_{s\rightarrow\infty}\overline{W}(s,t)\geq\lim_{s\rightarrow\infty}\underline{W}(s,t)$ for all $t>0$. Then $\overline{W}\geq\underline{W}$ in $[0,\infty)\times[0,\infty)$.
\end{lem}
\begin{rem}\label{remcomp}
The comparison principle above also holds for bounded space-time cylinders $[0,s_0]\times[0,t_0]$, if we assume that $\overline{W}$ and $\underline{W}$ belong to $C^0([0,s_0]\times[0,t_0])\cap C^{2,1}((0,s_0)\times(0,t_0))$ and fulfill the inequalities $\underline{W}(0,t)\leq\overline{W}(0,t)$ for $t\in[0,t_0]$, $\underline{W}(s_0,t)\leq\overline{W}(s_0,t)$ for $t\in[0,t_0]$ and $\underline{W}_0\leq\overline{W}_0$ in $[0,s_0]$ instead of the inequalities stated in the lemma above.
\end{rem}
In addition to allowing for this comparison principle, the approximate problem is also uniquely solvable in the classical sense.

\begin{lem}\label{solvWeps}
For $1>R>0,\rho\in(0,\frac{R}{2}), n>\alpha>2$ and $f_0>0$ let $F(s)$ be defined as in \eqref{Fdef}. Assume that $W_0$ satisfies \eqref{w0prop} with some $\mu>0$. Then for each $\varepsilon\in(0,1)$ there exists a unique function $W^{\ep}\in C^0([0,\infty)\times[0,\infty))\cap C^{2,1}((0,\infty)\times(0,\infty))$ which satisfies \eqref{WDGLeps} in the classical sense.
\end{lem}
Taking the limit $\varepsilon\searrow0$ to obtain $W^{\ep}\nearrow W$ in $(0,\infty)\times(0,\infty)$, which -- adopting the notion of \cite{TW11} and \cite{gv97} -- we will call the proper solution of \eqref{WDGL}. This limit procedure combined with the backwards transformation in \eqref{backtrafo} then results in a desired radial weak solution in the sense of Definition \ref{weakrad}.
\begin{lem}
For $1>R>0,\rho\in(0,\frac{R}{2}), n>\alpha>2$ and $f_0>0$ let $F(s)$ be defined as in \eqref{Fdef}. Assume that $u_0\in\Lo^1(\R^n)\cap\Lo^\infty(\mathbb{R}^n)$ is a radially symmetric function satisfying $\mu:=\sintRn u_0\intd x<\infty$. Then there exists at least one radial weak solution $u$ of \eqref{KS0f} in the sense of Definition \ref{weakrad}. Such a solution can be obtained by applying the backwards transformation \eqref{backtrafo} to the proper solution of \eqref{WDGL}.
\end{lem}
\setcounter{gleichung}{0} 
\section{Immediate blow-up of radial weak solutions}\label{cmeasvals2}
We require two further preparatory results not connected to the solution. The first is a variation of Gronwall's Lemma and can also be found in \cite[Lemma 2.1]{TW11}, whereto we refer for proof once again.
\begin{lem}\label{grönvar}
Suppose that $\Phi\in\W^{1,\infty}_{loc}(\mathbb{R})$ is non-decreasing, and that for some $t_1\in\mathbb{R},T>0$ and $c\in\mathbb{R}$ we are given two functions $y\in C^0([t_1,t_1+T))$ and $z\in C^1([t_1,t_1+T))$ such that
\begin{align*}
y(t)\geq c+\int\limits_{t_1}^t\Phi(y(\tau))\intd\tau\quad\text{for all }t\in(t_1,t_1+T)
\end{align*}
and
\begin{align*}
\left\{\begin{array}{c@{\,}l@{\quad}l@{\,}c}
z'(t)&=\Phi(z(t))\ &\text{for all }t\in(t_1,t_1+T),\\
z(t_1)&=c.&
	\end{array}\right.
\end{align*}
Then
\begin{align*}
y(t)\geq z(t)\quad\text{for all }t\in(t_1,t_1+T).
\end{align*}
\end{lem} 
The next lemma will provide us with functions fulfilling the role of test functions later in the proof of our main theorem. This lemma is an adjusted version of the corresponding construction from \cite[Lemma 2.2]{TW11}.
\begin{lem}\label{lemtestfunc}
For $1>R>0,\rho\in(0,\frac{R}{2}), n>\alpha>2$ and $f_0>0$ satisfying
\begin{align}\label{fproperty}
\owncount
f_0>\frac{2n}{\alpha}(n-2)(n-\alpha),
\end{align}
let $F(s)$ be defined as in \eqref{Fdef}. Furthermore, we define $h(n,\alpha,f_0):=(n-\alpha)(3n-4)-f_0$ and fix $\xi\in(4-\frac{4}{n},4]$ as well as $\delta\in(0,1)$ fulfilling the condition
\begin{align}\label{dproperty}
\owncount
\delta>\max\left\lbrace\frac{n-\alpha}{n},\frac{h(n,\alpha,f_0)+\sqrt{h(n,\alpha,f_0)^2+4f_0(n-\alpha)^2}}{2n(n-\alpha)}\right\rbrace.
\end{align}
Then there exist positive constants $a,b,k_0$ and $K_0$ depending only on $n,f_0,\alpha$ and $\xi$, such that for any $\gamma>\frac{4}{R-\rho}$ the function $\varphi^{(\gamma)}=:\varphi:(0,\infty)\rightarrow\R$ defined by
\begin{align}\label{phidef}
\owncount
\varphi(s):=\begin{cases}
\frac{a}{\gamma^\delta}s^{-\delta}-b&\text{ if }s< \frac{\xi}{\gamma},\\
e^{-\gamma s}&\text{ if }s\geq \frac{\xi}{\gamma},
\end{cases}
\end{align}
is positive, belongs to $\W^{2,\infty}_{loc}((0,\infty))$ and satisfies
\begin{align}\label{phiodi}
\owncount
n^2s^{\frac{2n-2}{n}}\varphi_{ss}+4(n^2-n)s^{\frac{n-2}{n}}\varphi_s-nF\varphi_s-nF_s\varphi\geq k_0\gamma^{\frac{2}{n}}\varphi\ \text{ a.e in }(0,\infty),
\end{align}
as well as
\begin{align}\label{phiint}
\owncount
\int\limits_0^{\infty}\frac{\varphi^2(s)}{\vert\varphi_s(s)\vert}\intd s\leq\frac{K_0}{\gamma^2}.
\end{align}
\end{lem}
\begin{bew}
Let us first verify that the assumed property \eqref{fproperty} for $f_0$ ensures that $\delta\in(0,1)$ can be chosen to satisfy \eqref{dproperty}. Clearly $\frac{n-\alpha}{n}<1$, thus we only have to verify that \eqref{fproperty} implies
\begin{align*}
h(n,\alpha,f_0)+\sqrt{h(n,\alpha,f_0)^2+4f_0(n-\alpha)^2}<2n(n-\alpha).
\end{align*}
Multiplication of \eqref{fproperty} with $\frac{4\alpha}{n-\alpha}$ implies
\begin{align*}
8n(n-2)<\frac{4\alpha}{n-\alpha}f_0.
\end{align*}
Adding $n^2+\frac{f_0^2}{(n-\alpha)^2}+\frac{8-4\alpha-2n}{n-\alpha}f_0-8n+16$ to both sides we obtain
\begin{align*}
9n^2+\frac{f_0^2}{(n-\alpha)^2}-24n+\frac{8-4\alpha-2n}{n-\alpha}f_0+16< n^2+\frac{f_0^2}{(n-\alpha)^2}-8n+\frac{8-2n}{n-\alpha}f_0+16,
\end{align*}
which may also be expressed as
\begin{align*}
\left(3n-4-\frac{f_0}{n-\alpha}\right)^2+4f_0 <\left(4-n+\frac{f_0}{n-\alpha}\right)^2.
\end{align*}
This implies
\begin{align*}
\sqrt{\left(3n-4-\frac{f_0}{n-\alpha}\right)^2+4f_0}<4-n+\frac{f_0}{n-\alpha}
\end{align*}
and thus, multiplying with $(n-\alpha)$, recalling the definition of $h(n,\alpha,f_0)$ and reordering the terms appropriately, we obtain 
\begin{align*}
h(n,\alpha,f_0)+\sqrt{h(n,\alpha,f_0)^2+4f_0(n-\alpha)^2}<2n(n-\alpha).
\end{align*}
Next, we observe that the fact $\gamma>\frac{4}{R-\rho}$ implies that $\xi$ satisfies the inequalities $(R-\rho)\gamma>\xi>4-\frac{4}{n}>\delta$. Setting $a:=\frac{\xi^{\delta+1}}{\delta}e^{-\xi}$ and $b:=\left(\frac{\xi}{\delta}-1\right)e^{-\xi}$ we have $a>0$ and $b>0$. With the constants $a,b$ defined this way, the function $\varphi(s)$ in \eqref{phidef} fulfills
\begin{align*}
\lim_{s\nearrow\frac{\xi}{\gamma}}\varphi(s)=\frac{a}{\gamma^\delta}\frac{\gamma^\delta}{\xi^\delta}-b=\frac{\xi}{\delta}e^{-\xi}-\left(\frac{\xi}{\delta}-1\right)e^{-\xi}=e^{-\xi}=\varphi\left(\frac{\xi}{\gamma}\right)
\end{align*}
and thus, $\varphi$ is continuous in $(0,\infty)$. Similarly, simple calculation shows
\begin{align*}
\varphi_s(s)=\begin{cases}
\frac{-a\delta}{\gamma^\delta}s^{-(\delta+1)}&\text{ if }s< \frac{\xi}{\gamma},\\
-\gamma e^{-\gamma s}&\text{ if }s> \frac{\xi}{\gamma}.
\end{cases}
\end{align*}
Therefore, $\varphi$ is monotonically decreasing in $(0,\infty)$, which by the definition of $\varphi$ immediately implies $\varphi>0$ in $(0,\infty)$. Moreover, it holds that
\begin{align*}
\lim_{s\nearrow\frac{\xi}{\gamma}}\varphi_s(s)=\frac{-a\delta}{\gamma^\delta}\frac{\gamma^{\delta+1}}{\xi^{\delta+1}}=-\gamma e^{-\xi}=\varphi_s\left(\frac{\xi}{\gamma}\right)
\end{align*}
and hence $\varphi_s$ is continuous in $(0,\infty)$ as well. One further differentiation shows
\begin{align*}
\varphi_{ss}(s)=\begin{cases}
\frac{a\delta(\delta+1)}{\gamma^\delta}s^{-(\delta+2)}&\text{ if }s< \frac{\xi}{\gamma},\\
\gamma^2 e^{-\gamma s}&\text{ if }s> \frac{\xi}{\gamma},
\end{cases}
\end{align*}
therefore clearly $\varphi\in\W^{2,\infty}_{loc}((0,\infty))$. Utilizing the expressions above, we see that for $s>\frac{\xi}{\gamma}$ holds:
\begin{align}\label{Lphiest}
\owncount
L\varphi(s):=\ &n^2s^{\frac{2n-2}{n}}\varphi_{ss}(s)+4(n^2-n)s^{\frac{n-2}{n}}\varphi_s(s)-nF(s)\varphi_s(s)-nF_s(s)\varphi(s)\nonumber\\
=\ &\left(n^2s^{\frac{2n-2}{n}}\gamma^2-4(n^2-n)s^{\frac{n-2}{n}}\gamma\right)\!e^{-\gamma s} +\left(n\gamma F(s)-nF_s(s)\right)\!e^{-\gamma s}
\end{align}
To estimate $n\gamma F(s)-nF_s(s)$ from below we distinguish the cases $s>R-\rho$ and $R-\rho>s>\frac{\xi}{\gamma}$. In the first case we have
\begin{align*}
n\gamma F(s)-nF_s(s)\geq n\gamma F(R-\rho)-nF_s(R-\rho)=\left(\frac{n\gamma}{n-\alpha}(R-\rho)-1\right)f_0(R-\rho)^{-\frac{\alpha}{n}},
\end{align*}
by the monotonicity of $F$ and $F_s$ and their explicit values at $s=R-\rho$ shown in \eqref{Fdef} and \eqref{Fsdef}, respectively. Because of $n>\alpha$ and $(R-\rho)\gamma>4-\frac{4}{n}>1$, we have $\frac{n\gamma(R-\rho)}{n-\alpha}>1$ and since $f_0>0$ and $R-\rho>0$ also hold, this shows 
\begin{align*}
n\gamma F(s)-nF_s(s)\geq0.
\end{align*}
In the case $R-\rho>s>\frac{\xi}{\gamma}$ we make use of $s^{\frac{n-\alpha}{n}}\geq\frac{\xi}{\gamma}s^{-\frac{\alpha}{n}}$, to obtain
\begin{align*}
n\gamma F(s)-nF_s(s)&=\left(\frac{n\gamma}{n-\alpha}s^{\frac{n-\alpha}{n}}-s^{-\frac{\alpha}{n}}\right)f_0\geq\left(\frac{n\xi}{n-\alpha}-1\right)f_0 s^{-\frac{\alpha}{n}}.
\end{align*}
This term is nonnegative as well, since $\frac{n-\alpha}{n}<1<\xi$ implies $\frac{n\xi}{n-\alpha}>1$.
Thus we may drop the last term in \eqref{Lphiest} and acquire
\begin{align*}
L\varphi(s)\geq\ &\left(n^2 s^{\frac{2n-2}{n}}\gamma^2-4(n^2-n)\gamma s^{\frac{n-2}{n}}\right)e^{-\gamma s}\\
=\ &\left(n^2s^{\frac{2n-2}{n}}\gamma^2-4(n^2-n)\gamma s^{\frac{n-2}{n}}\right)\!\varphi(s)\ \text{ for all }s>R-\rho.
\end{align*}
Using $s>\frac{\xi}{\gamma}$ again, we estimate $s^{\frac{2n-2}{n}}\geq\frac{\xi}{\gamma}s^{\frac{n-2}{n}}$, which leads to
\begin{align*}
L\varphi(s)\geq\left(n^2\xi-4(n^2-n)\right)\gamma s^{\frac{n-2}{n}}\!\varphi(s)\ \text{ for all }s>R-\rho.
\end{align*}
Here the choice of $\xi>4-\frac{4}{n}$ implies that the foremost factor is nonnegative and thus, using $s>\frac{\xi}{\gamma}$ once more, we obtain
\begin{align*}
L\varphi(s)\geq\ &\left(n^2\xi-4(n^2-n)\right)\gamma^{\frac{2}{n}}\xi^{\frac{n-2}{n}}\varphi(s)=:c_1\gamma^{\frac{2}{n}}\varphi(s)\ \text{ for all }s>R-\rho.
\end{align*}
On the other hand for $0<s<\frac{\xi}{\gamma}<R-\rho< 1$:
\begin{align*}
L\varphi(s)=\ &n^2\frac{a\delta(\delta+1)}{\gamma^\delta}s^{-\delta-\frac{2}{n}}-4(n^2-n)\frac{a\delta}{\gamma^\delta}s^{-\delta-\frac{2}{n}}+\frac{nf_0}{n-\alpha}\frac{a\delta}{\gamma^\delta}s^{-\delta-\frac{\alpha}{n}}-f_0s^{-\frac{\alpha}{n}}\left(\frac{a}{\gamma^\delta}s^{-\delta}-b\right)\\
\geq\ &\left(\left(n^2\delta(\delta+1)-4\delta(n^2-n)\right)s^{-\frac{2}{n}}+\left(\frac{\delta n}{n-\alpha}-1\right)f_0s^{-\frac{\alpha}{n}}\right)\frac{a}{\gamma^\delta}s^{-\delta}.
\intertext{Recalling the first argument of the maximum in \eqref{dproperty}, we have $\frac{\delta n}{n-\alpha}-1>0$ and since $1>s$ and $\alpha>2$ imply $s^{-\frac{\alpha}{n}}>s^{-\frac{2}{n}}$, we obtain}
L\varphi(s)\geq\ &\left(\left(n^2\delta(\delta+1)-4\delta(n^2-n)\right)s^{-\frac{2}{n}}+\left(\frac{\delta n}{n-\alpha}-1\right)f_0s^{-\frac{2}{n}}\right)\frac{a}{\gamma^\delta}s^{-\delta}\\
=\ &\left(n^2\delta^2+\left(\frac{n f_0}{(n-\alpha)}-3n^2+4n\right)\delta-f_0\right)s^{-\frac{2}{n}}\frac{a}{\gamma^\delta}s^{-\delta}\ \text{ for all }0<s<\frac{\xi}{\gamma}.
\end{align*}
Observing that the larger root of the equation $n^2\delta^2+\left(\frac{n f_0}{(n-\alpha)}-3n^2+4n\right)\delta-f_0=0$ for $\delta$ equals the second argument of the maximum in \eqref{dproperty}, the choice of $\delta$ implies that $n^2\delta^2+\left(\frac{n f_0}{(n-\alpha)}-3n^2+4n\right)\delta-f_0>0$ holds. Thus the factor in front is positive and we may estimate
\begin{align*}
L\varphi(s)\geq\ &\left(n^2\delta^2+\left(\frac{n f_0}{(n-\alpha)}-3n^2+4n\right)\delta-f_0\right)s^{-\frac{2}{n}}\varphi\ \text{ for all }0<s<\frac{\xi}{\gamma}.
\intertext{Now we can use $s<\frac{\xi}{\gamma}$ to obtain}
L\varphi(s)\geq\ &\left(n^2\delta^2+\left(\frac{n f_0}{(n-\alpha)}-3n^2+4n\right)\delta-f_0\right)\xi^{-\frac{2}{n}}\gamma^{\frac{2}{n}}\varphi=:c_2\gamma^{\frac{2}{n}}\varphi\ \text{ for all }0<s<\frac{\xi}{\gamma}.
\end{align*}
Choosing $k_0:=\min\left\lbrace c_1,c_2\right\rbrace$, the asserted inequality \eqref{phiodi} holds a.e in $(0,\infty)$. To show \eqref{phiint} we use $b>0$ again to calculate: 
\begin{align*}
\int\limits_0^{\frac{\xi}{\gamma}}\frac{\varphi^2}{\vert\varphi_s\vert}\intd s&\leq\int\limits_0^{\frac{\xi}{\gamma}}\frac{\frac{a^2}{\gamma^{2\delta}}s^{-2\delta}}{\frac{a\delta}{\gamma^\delta}s^{-\delta-1}}\intd s=\int\limits_0^{\frac{\xi}{\gamma}}\frac{a}{\delta\gamma^\delta}s^{1-\delta}\intd s=\frac{a\xi^{2-\delta}}{\delta(2-\delta)\gamma^2}
\end{align*}
and
\begin{align*}
\int\limits_{\frac{\xi}{\gamma}}^{\infty}\frac{\varphi^2}{\vert\varphi_s\vert}\intd s&=\int\limits_{\frac{\xi}{\gamma}}^{\infty}\frac{e^{-2\gamma s}}{\gamma e^{-\gamma s}}\intd s=\frac{e^{-\xi}}{\gamma^2}.
\end{align*}
Hence the asserted statement in \eqref{phiint} holds for $K_0:=\left(\frac{a\xi^{2-\delta}}{2-\delta}+e^{-\xi}\right)$, which completes the proof.
\end{bew}

Before we begin with the proof of our main result let us fix some parameters.

\begin{lem}\label{lempara}
Assume that the conditions of Lemma \ref{lemtestfunc} hold and $W_0$ satisfies \eqref{w0prop} for some $\mu>0$, as well as $W_0\in C^0([0,1])\cap C^2((0,1))$, $W_0(s)\geq c_0s$ for some $c_0>0$ and all $s\in[0,1]$, $W_{0s}(s)\geq0$ for all $s\in(0,1)$ and $W_{0ss}(s)\geq0$ for all $s\in(0,1)$. Denote by $W$ the corresponding proper solution of \eqref{WDGL}. Furthermore, we set $\xi=4$ and let $a,b,k_0,K_0$ be the positive constants according to Lemma \ref{lemtestfunc}. Then for any $t_0\geq0$ and $\eta>0$ we can find $\kappa>0$ fulfilling
\begin{align}\label{kappaineq}
\owncount
\kappa\leq\frac{k_0\eta}{8}
\end{align}
and  $\gamma>\frac{4}{R-\rho}$ such that
\begin{align}\label{gammaineqe}
\owncount
1+\frac{2k_0K_0e^{\kappa\gamma^{\frac{2}{n}}}}{W(\kappa\gamma^{\frac{2-n}{n}},t_0+\frac{\eta}{2})\gamma^{\frac{n-2}{n}}}\leq e^{2\kappa\gamma^{\frac{2}{n}}}
\end{align}
as well as
\begin{align}\label{gammaineqR}
\owncount
\gamma>\left(\frac{\xi}{\kappa}\right)^{\frac{n}{2}}
\end{align}
hold.
\end{lem}
\begin{bew}
Because of $\eta>0$ and $k_0>0$, we may choose a sufficiently small $\kappa>0$ which satisfies the inequality \eqref{kappaineq}. For $\varepsilon\in(0,1)$ let $W^{\ep}$ denote the solution to \eqref{WDGLeps}. We want to estimate $W^{\ep}$ from below in $[0,1]\times[0,t_0+\frac{\eta}{2}]$ by a suitable subsolution. To this end we observe that we have $W_0(1)>0$, since $W_0(s)\geq c_0 s$ for every $s\in[0,1]$. Thus, by the strong maximum principle applied to \eqref{WDGLeps}, $W^{\ep}$ is positive in $(\frac{1}{4},\infty)\times(0,\infty)$. Hence, the number 
\begin{align*}
c_1:=\inf_{\tau\in[0,t_0+\frac{\eta}{2}]}\left\{1,\frac{W^{\ep}(\frac{1}{2},\tau)}{W_0(1)}\right\}
\end{align*}
is positive and well-defined. In fact, setting $\Wu(s,t):=c_1 s^2 W_0(s)$ we see that
\begin{align*}
\Wu(0,t)&=0=W^{\ep}(0,t)\text{ for all }t\in[0,t_0+\frac{\eta}{2}],\\
\Wu(1,t)&=c_1 W_0(1)\leq W^{\ep}(\tfrac{1}{2},t)\leq W^{\ep}(1,t)\text{ for all }t\in[0,t_0+\frac{\eta}{2}],
\intertext{and}
\Wu(s,0)&=c_1 s^2W_0(s)\leq W_0(s),\text{ for }s\in[0,1].
\end{align*}
Furthermore we have
\begin{align*}
\mbox{ }&n^2s^{\frac{2n-2}{n}}\Wu_{ss}+\chie\Wu\,\Wu_s+n\chie F\Wu_s\\=&\ n^2s^{\frac{2n-2}{n}}\left(2c_1 W_0+c_1 s^2W_{0ss}+4sc_1 W_{0s}\right)+\left(c_1 s^2W_0+nF\right)\chie\left(2c_1 sW_0+c_1 s^2W_{0s}\right)\\
\geq&\ 0=\Wu_{\,t}
\end{align*}
in $(0,1)\times(0,t_0+\tfrac{\eta}{2})$, since $W_{0s}(s)\geq0$ and $W_{0ss}(s)\geq0$ for all $s\in(0,1)$. Thus, we may use the comparison principle, see Remark \ref{remcomp}, to deduce
\begin{align*}
W^{\ep}(s,t)\geq\Wu(s,t)\text{ in }[0,1]\times[0,t_0+\frac{\eta}{2}].
\end{align*}
In particular, since $W^{\ep}\nearrow W$ as $\varepsilon\searrow0$, we have
\begin{align}\label{Wsubineq}
\owncount
W(s,t_0+\frac{\eta}{2})\geq W^{\ep}(s,t_0+\frac{\eta}{2})\geq c_1 s^2 W_0(s)\geq c_1c_0 s^3=:p(s)\text{ for }s\in[0,1]
\end{align}
since $W_0(s)\geq c_0 s$ in $[0,1]$. This inequality at hand we can now verify that by choosing $\gamma$ sufficiently large \eqref{gammaineqe} is indeed fulfilled. To this end, let $s_0\in(0,1)$ be so small, such that the inequalities
\begin{align}\label{kappa1}
\owncount
\frac{k_0K_0}{\kappa}\leq c_0c_1s_0^3\sinh\left(\kappa\left(\frac{\kappa}{s_0}\right)^{\frac{2}{n-2}}\right)
\end{align}
and
\begin{align}\label{kappa2}\owncount
s_0<\left(\frac{2\kappa^{\frac{n}{n-2}}}{3(n-2)}\right)^{\frac{n-2}{2}}
\end{align}
hold. Rearranging \eqref{kappa2} we see that
\begin{align*}
\frac{3}{2}s^2<\frac{\kappa^{\frac{n}{n-2}}}{n-2}s^{\frac{2n-6}{n-2}}
\end{align*}
holds for all $s\in(0,s_0)$. Hence,
\begin{align*}
\frac{3}{2}s^2 \exp\left(\kappa\left(\frac{\kappa}{s}\right)^{\frac{2}{n-2}}\right)<\frac{\kappa^{\frac{n}{n-2}}}{n-2}s^{\frac{2n-6}{n-2}}\exp\left(\kappa\left(\frac{\kappa}{s}\right)^{\frac{2}{n-2}}\right)
\end{align*}
for all $s\in(0,s_0)$. This implies that $g(s):=s^3\sinh\left(\frac{\kappa^{\frac{n}{n-2}}}{s^{\frac{2}{n-2}}}\right)$ is monotonically decreasing, since
\begin{align*}
\frac{\partial g(s)}{\partial s}=
\left(\frac{3}{2}s^2-\frac{\kappa^{\frac{n}{n-2}}}{n-2}s^{\frac{2n-6}{n-2}}\right)\exp\left(\kappa\left(\frac{\kappa}{s}\right)^{\frac{2}{n-2}}\right)\!-\!\left(\frac{3}{2}s^2+\frac{\kappa^{\frac{n}{n-2}}}{n-2}s^{\frac{2n-6}{n-2}}\right)\exp\left(-\kappa\left(\frac{\kappa}{s}\right)^{\frac{2}{n-2}}\right)
\leq0
\end{align*}
for all $s\in(0,s_0)$. Making use of this monotonicity and \eqref{kappa1} we obtain
\begin{align*}
c_0 c_1 s^3\sinh\left(\kappa\left(\frac{\kappa}{s}\right)^{\frac{2}{n-2}}\right)\geq c_0 c_1 s_0^3\sinh\left(\kappa\left(\frac{\kappa}{s_0}\right)^{\frac{2}{n-2}}\right)\geq\frac{k_0 K_0}{\kappa}\text{ for all }s\in(0,s_0).
\end{align*}
Utilizing $s_0<1$ this yields
\begin{align}\label{polyineq}
\owncount
p(s)\geq\frac{k_0K_0s}{\kappa\sinh\left(\kappa\left(\frac{\kappa}{s}\right)^{\frac{2}{n-2}}\right)}
\end{align}
for all $s\in(0,s_0)$. Hence, taking $\gamma>\frac{4}{R-\rho}$ such that $s_0>\kappa\gamma^{\frac{2-n}{n}}$ and $\gamma>\left(\tfrac{\xi}{\kappa}\right)^{\frac{n}{2}}$ are fulfilled, we obtain \eqref{gammaineqR} as well as
\begin{align*}
W\left(\kappa\gamma^{\frac{2-n}{n}},t_0+\frac{\eta}{2}\right)\geq p(\kappa\gamma^{\frac{2-n}{n}})\geq\frac{2k_0K_0}{\gamma^{\frac{n-2}{n}}}\frac{e^{\kappa\gamma^{\frac{2}{n}}}}{e^{2\kappa\gamma^{\frac{2}{n}}}-1},
\end{align*}
from \eqref{Wsubineq} and \eqref{polyineq}. \eqref{gammaineqe} then follows from reordering the inequality above.
\end{bew}

Inspired by the methods of \cite[Lemma 2.3]{TW11}, we can now use functions of the type described in Lemma \ref{lemtestfunc} to prove that the spatial derivative of the proper solution $W$ of \eqref{WDGL} blows up immediately:
\begin{prop}\label{wspropo}
Suppose that the conditions of Lemma \ref{lemtestfunc} hold. Then for any $\beta\geq1$ and each $t_0\geq0$ the proper solution $W$ of \eqref{WDGL} fulfills
\begin{align}\label{wsblow}
\owncount
\sup_{\substack{s\,>\,0\\t\in\,(t_0,t_0+\eta)}}\frac{W(s,t)}{s^\beta}=\infty\quad\text{for all }\eta>0.
\end{align}
In particular, we have
\begin{align}\label{wsblow2}
\owncount
\|W_s\|_{\Lo^\infty((0,\infty)\times(t_0,t_0+\eta))}=\infty\quad\text{for all }\eta>0.
\end{align}
\end{prop}
\begin{bew}
We work along the lines of a contradiction argument and assume to this end that there exist $\beta\geq1, t_0\geq0$, $\eta>0$ and $c>0$ such that
\begin{align}\label{Wscontra}
\owncount
W(s,t)\leq c s^\beta\quad\text{for all }s\geq0\text{ and }t\in[t_0,t_0+\eta].
\end{align}
Since $\beta\geq1$, we can fix $\delta\in(0,1)$ satisfying $\beta>\delta$ as well as the property \eqref{dproperty}. Furthermore set $\xi=4$ and let $a,b,k_0$ and $K_0$ be the positive constants defined in Lemma \ref{lemtestfunc}. Corresponding to these parameters let $\kappa,\gamma$ be the positive constants given by Lemma \ref{lempara}.

With these parameters we define $\varphi:=\varphi^{(\gamma)}$ as in \eqref{phidef} of Lemma \ref{lemtestfunc}. In particular, $\varphi$ fulfills the differential inequality \eqref{phiodi} and the integral inequality \eqref{phiint}. Both will be required later in this proof.  Recalling the cut-off functions $\chie$ for $\varepsilon\in(0,1)$, mentioned in \eqref{cut-off1} and \eqref{cut-off2}, we multiply the approximation problem \eqref{WDGLeps} by $\chie\varphi$. Integration by parts over $s\in(0,s_0)$, where $s_0$ is an arbitrary number satisfying $s_0>\max\left\lbrace\frac{\xi}{\gamma},\frac{2n-2}{n\gamma},\varepsilon\right\rbrace$, results in
\begin{align}\label{blowupproofeq1}
\owncount
\frac{\intd}{\intd t}\intos\chie\varphi W^{\ep}\intd s&=\intos\chie\varphi\left[n^2 s^{\frac{2n-2}{n}}W_{ss}^{\ep}+\frac{1}{2}\chie\left((\W^{\ep})^2\right)_s+n\chie F\W_s^{\ep}\right]\intd s\nonumber\\
&=\ n^2\intos\left(s^{\frac{2n-2}{n}}\chie\varphi\right)_{ss} W^{\ep}\intd s-\frac{1}{2}\intos\left((\chie)^2\varphi\right)_s(W^{\ep})^2\intd s\nonumber\\
&\mbox{ }\quad-n\intos\left((\chie)^2\varphi F\right)_s W^{\ep}\intd s+B(t)\text{ for }t>0,
\end{align}
where $B(t)$ are the collected boundary terms, that is
\begin{align}\label{blowupproofeq2}
\owncount
B(t):=\bigg[n^2 s^{\frac{2n-2}{n}}\chie\varphi W^{\ep}_s&-n^2\left(s^{\frac{2n-2}{n}}\chie\varphi\right)_s W^{\ep}\nonumber\\
\hspace*{80pt}&
+\frac{1}{2}(\chie)^2(W^{\ep})^2\varphi+n(\chie)^2  W^{\ep}\varphi F\bigg]_0^{s_0}.
\end{align}
Calculating the mixed derivative term
\begin{align}\label{blowupproofeq3}
\owncount
n^2\left( s^{\frac{2n-2}{n}}\varphi\chie\right)_s&=n^2s^{\frac{2n-2}{n}}\varphi_s\chie+n^2s^{\frac{2n-2}{n}}\varphi\chi^{\ep}_s+n(2n-2)s^{\frac{n-2}{n}}\varphi\chie,
\end{align}
as well as its companions
\begin{align}
\owncount\label{blowupproofeq5}
\frac{1}{2}\left(\varphi(\chie)^2\right)_s&=\frac{1}{2}(\chie)^2\varphi_s+\varphi\chie\chi^{\ep}_s, \\
\owncount\label{blowupproofeq6}
n\left((\chie)^2\varphi F\right)_s&=n(\chie)^2\varphi_s F+2n\chie\chi^{\ep}_s\varphi F+n(\chie)^2\varphi F_s,
\intertext{and}
\owncount\label{blowupproofeq4}
n^2\left( s^{\frac{2n-2}{n}}\varphi\chie\right)_{ss}&=n^2s^{\frac{2n-2}{n}}\varphi_{ss}\chie+2n^2s^{\frac{2n-2}{n}}\varphi_s\chi^{\ep}_s+4n(n-1)s^{\frac{n-2}{n}}\varphi_s\chie\nonumber\\
&\quad\!+n^2s^{\frac{2n-2}{n}}\varphi\chi^{\ep}_{ss}+4n(n-1)s^{\frac{n-2}{n}}\varphi\chi^{\ep}_s+(2n^2-6n+4)s^{-\frac{2}{n}}\varphi\chie,
\end{align}
we can make use of the facts $\chie\equiv0$ on $[0,\frac{\varepsilon}{2}]$, $\chie\equiv1$ on $[\varepsilon,\infty)$, $\varphi(s)=e^{-\gamma s}$ on $[\frac{\xi}{\gamma},\infty)$ and $s_0>\max\lbrace\frac{2n-2}{n\gamma},\frac{\xi}{\gamma},\varepsilon\rbrace$, to express \eqref{blowupproofeq2} as
\begin{align*}
B(t)=n^2 s_0^{\frac{2n-2}{n}}e^{-\gamma s_0}W^{\ep}_s\!(s_0,t)&+\left(n^2s_0^{\frac{2n-2}{n}}\gamma e^{-\gamma s_0}-n(2n-2)s_0^{\frac{n-2}{n}}e^{-\gamma s_0}\right)W^{\ep}\!(s_0,t)\\
&+\frac{1}{2} (W^{\ep}\!(s_0))^2e^{-\gamma s_0}+n W^{\ep}\!(s_0,t)e^{-\gamma s_0}F(s_0)\ \text{ for all }t>0.
\end{align*}
Recalling that $W^{\ep}_s\geq0$ and $F\geq0$ we can drop multiple nonnegative terms to obtain
\begin{align*}
B(t)\geq\left(n\gamma s_0-(2n-2)\right)ns_0^{\frac{n-2}{n}} e^{-\gamma s_0}W^{\ep}\!(s_0,t)\ \text{ for all }t>0,
\end{align*}
which by choice of $s_0$ then implies $B(t)\geq0\ \text{ for all }t>0$. Thus, inserting \eqref{blowupproofeq5}--\eqref{blowupproofeq4} into \eqref{blowupproofeq1}, we get
\begin{align}\label{blowupproofeq7}
\owncount
\frac{\intd}{\intd t}\intos\varphi\chie W^{\ep}\intd s&\geq\intos\chie W^{\ep}\bigg[n^2 s^{\frac{2n-2}{n}}\varphi_{ss}+4n(n-1)s^{\frac{n-2}{n}}\varphi_s\nonumber\\
&\hspace*{55pt}+(2n^2-6n+4)s^{-\frac{2}{n}}\varphi-n\chie\varphi_s F-n\chie\varphi F_s\bigg]\intd s\nonumber\\
&-\frac{1}{2}\intos\varphi_s(\chie)^2(W^{\ep})^2\intd s+I_1(t)+I_2(t)\ \text{ for all }t>0,
\end{align}
where we set $I_1(t):=\intos\chi^{\ep}_sW^{\ep}\left(2n^2s^{\frac{2n-2}{n}}\varphi_s+4n(n-1)s^{\frac{n-2}{n}}\varphi-\varphi\chie W^{\ep}-2n\chie \varphi F\right)\intd s$ and $I_2(t):=\intos n^2 s^{\frac{2n-2}{n}}\varphi\chi^{\ep}_{ss}W^{\ep}\intd s$. Using the properties of the cut-off function $\chie$ we can estimate both $I_1$ and $I_2$ from below.
For that, we first recall that $\chi^{\ep}_s\leq\frac{c_\chi}{\varepsilon},
\chie\leq1,\chi^{\ep}_s\equiv0$ on $(\varepsilon,\infty)$ and $W^{\ep}\leq W\leq\frac{n\mu}{\vert S_{n-1}\vert}\ \text{ for all }s\geq0$ and $t>0$. Thus, using $F(s)\leq F(R+\rho)$ for all $s>0$, we can estimate
\begin{align}
\owncount\label{blowupproofeq8}
I_1(t)\ &\geq\intos\left(2n^2s^{\frac{2n-2}{n}}\varphi_s-\varphi\chie W^{\ep}-2n\chie \varphi F\right)\chi^{\ep}_s W^{\ep}\intd s\nonumber\\
&\geq-\int\limits_0^{\varepsilon}\left(2n^2s^{\frac{2n-2}{n}}\vert\varphi_s\vert+\varphi W+2n\varphi F(R+\rho)\right)W\frac{c_\chi}{\varepsilon}\intd s.\nonumber\\
\intertext{Next, as long as $\varepsilon<\frac{\xi}{\gamma}$, we make use of the definitions of $\varphi, F$ and our contradiction assumption $W\leq c s^\beta$ for all $t\in(t_0,t_0+ \eta)$ and $s\geq0$ in \eqref{Wscontra}, to obtain}
I_1(t)\ &\geq-\int\limits_0^{\varepsilon}\left(2n^2\frac{a\delta}{\gamma^\delta}s^{\frac{2n-2}{n}-\delta-1}+\frac{a}{\gamma^\delta}s^{-\delta}\frac{n\mu}{\vert S_{n-1}\vert}+\frac{2nf_0}{n-\alpha}(R+\rho)^{\frac{n-\alpha}{n}}\frac{a}{\gamma^\delta}s^{-\delta}\right)\frac{c_\chi}{\varepsilon}cs^\beta\intd s.\nonumber\\
\intertext{Since $\varepsilon<1$ and $(R+\rho)^{\frac{n-\alpha}{n}}<2^{\frac{n-\alpha}{n}}<2$ one further estimation shows}
I_1(t)\ &\geq-\frac{cc_\chi a}{\gamma^\delta\varepsilon}\left(2n^2+\frac{n\mu}{\vert S_{n-1}\vert}+\frac{4nf_0}{n-\alpha}\right)\int\limits_0^{\varepsilon}s^{\beta-\delta}\intd s\ \text{ for all }t\in(t_0,t_0+\eta)
\end{align}
as long as $\varepsilon<\frac{\xi}{\gamma}$. Similarly, recalling $\vert\chi^{\ep}_{ss}\vert\leq\frac{c_\chi}{\varepsilon^2}$, we have
\begin{align}
\owncount\label{blowupproofeq9}
\vert I_2(t)\vert&\leq n^2\frac{cc_\chi}{\varepsilon^2}\frac{a}{\gamma^\delta}\int\limits_0^{\varepsilon} s^{\frac{2n-2}{n}+\beta-\delta}\intd s\nonumber\\
&\leq\frac{n^2c c_\chi a}{\gamma^\delta\varepsilon^2}\int\limits_0^\varepsilon s^{\beta-\delta+1}\intd s\ \text{ for all }t\in(t_0,t_0+\eta)
\end{align}
as long as $\varepsilon<\frac{\xi}{\gamma}$. By the choice of $\delta<\beta$, the combination of \eqref{blowupproofeq8} and \eqref{blowupproofeq9} yields $(I_1+I_2)(t)\geq-c_2\varepsilon^{\beta-\delta}$ in $(t_0,t_0+\eta)$ for some $c_2>0$ independent of $\varepsilon$, as long as $\varepsilon<\frac{\xi}{\gamma}$, so that \eqref{blowupproofeq7} reduces to
\begin{align*}
&\frac{\intd}{\intd t}\intos\varphi\chie W^{\ep}\intd s\\
\geq\ &\intos\chie W^{\ep}\bigg[n^2 s^{\frac{2n-2}{n}}\varphi_{ss}+4n(n-1)s^{\frac{n-2}{n}}\varphi_s+(2n^2-6n+4)s^{-\frac{2}{n}}\varphi\\
&\hspace*{70pt}-n\chie\varphi_s F-n\chie\varphi F_s\bigg]\intd s-\frac{1}{2}\intos\varphi_s(\chie)^2(W^{\ep})^2\intd s-c_2 \varepsilon^{\beta-\delta}\\
\geq\ &\intos\chie W^{\ep}\bigg[n^2 s^{\frac{2n-2}{n}}\varphi_{ss}+4n(n-1)s^{\frac{n-2}{n}}\varphi_s-nF\varphi_s\chie-nF_s\varphi\chie\bigg]\intd s\\
&\hspace*{70pt}-\frac{1}{2}\intos\varphi_s(\chie)^2(W^{\ep})^2\intd s-c_2 \varepsilon^{\beta-\delta}\ \text{ for all }t\in(t_0,t_0+\eta),
\end{align*}
as long as $\varepsilon<\frac{\xi}{\gamma}$. Now \eqref{phiodi} of Lemma \ref{lemtestfunc} implies
\begin{align*}
n^2 s^{\frac{2n-2}{n}}\varphi_{ss}\!+4n(n-1)s^{\frac{n-2}{n}}\varphi_s\!-\!(nF\varphi_s+nF_s\varphi)\chie\!\geq k_0\gamma^{\frac{2}{n}}\varphi+(1-\chie)(nF\varphi_s+nF_s\varphi)
\end{align*}
a.e in $(0,\infty)$ and thus we obtain
\begin{align*}
&\frac{\intd}{\intd t}\intos\varphi\chie W^{\ep}\intd s\\
\geq\ &k_0\gamma^{\frac{2}{n}}\intos\varphi\chie W^{\ep}\intd s+\!\intos\! n(1-\chie)F\varphi_s\chie W^{\ep}\intd s+\!\intos\! n(1-\chie)F_s\varphi\chie W^{\ep}\intd s\\
&-\frac{1}{2}\intos\varphi_s(\chie)^2(W^{\ep})^2\intd s-c_2\varepsilon^{\beta-\delta}\ \text{ for all }t\in(t_0,t_0+\eta),
\end{align*}
for $\varepsilon<\frac{\xi}{\gamma}$. Setting $y^{(\varepsilon,s_0)}(t):=\intos\varphi(s)\chie(s)W^{\ep}\!(s,t)\intd s$ and integrating over $(t_0+\frac{\eta}{2},t)=:(t_1,t)$, the inequality above takes the form
\begin{align*}
y^{(\varepsilon,s_0)}(t)&\geq y^{(\varepsilon,s_0)}(t_1)+k_0\gamma^{\frac{2}{n}}\int\limits_{t_1}^{t}\!\intos\varphi\chie W^{\ep}\intd s\intd t-\frac{1}{2}\int\limits_{t_1}^{t}\!\intos\varphi_s(\chie)^2(W^{\ep})^2\intd s\intd t\\
&+n\int\limits_{t_1}^{t}\!\intos F\varphi_s(1-\chie)\chie W^{\ep}\intd s\intd t+n\int\limits_{t_1}^{t}\!\intos F_s\varphi(1-\chie)\chie W^{\ep}\intd s\intd t\\
&-c_3\varepsilon^{\beta-\delta}(t-t_1)\ \text{ for all }t\in(t_1,t_0+\eta)\text{ and }0<\varepsilon<\min\left\lbrace1,\frac{\xi}{\gamma}\right\rbrace.
\end{align*}
Observing that $\vert\intos\varphi\chie W^{\ep}\intd s\vert\leq\left(\frac{a\xi^{1-\delta}}{(1-\delta)\gamma}+\frac{e^{-\xi}}{\gamma}\right)\frac{\mu n}{\vert S_{n-1}\vert}$ holds for all $s_0\in(0,\infty)$ we may use the monotone convergence theorem to take $s_0\nearrow\infty$ and obtain
\begin{align*}
\displaystyle\lim_{s_0\nearrow\infty}\int\limits_{t_1}^{t}\!\intos\varphi\chie W^{\ep}\intd s\intd t=\int\limits_{t_1}^{t}\!\int\limits_0^{\infty}\varphi\chie W^{\ep}\intd s\intd t<\infty\ \text{ for all }t\in(t_1,t_0+\eta).
\end{align*}
In a similar fashion -- using not only the exponential decay of $\varphi$ but also of $\varphi_s$ -- we can apply the monotone convergence theorem for
\begin{align*}
\displaystyle\lim_{s_0\nearrow\infty}&\int\limits_{t_1}^{t}\!\intos\varphi_s(\chie)^2(W^{\ep})^2\intd s\intd t,\\
\lim_{s_0\nearrow\infty}&\int\limits_{t_1}^{t}\!\intos F\varphi_s(1-\chie)\chie W^{\ep}\intd s\intd t
\intertext{and}
\displaystyle\lim_{s_0\nearrow\infty}&\int\limits_{t_1}^{t}\!\intos F_s\varphi(1-\chie)\chie W^{\ep}\intd s\intd t
\end{align*}
to see that the function $\displaystyle y^{\ep}\!(t):=\lim_{s_0\rightarrow\infty}y^{(\varepsilon,s_0)}(t)=\int\limits_0^{\infty}\varphi(s)\chie(s)W^{\ep}(s,t)\intd s$ satisfies the inequality
\begin{align}\label{blowupproofeq10}
\owncount
y^{\ep}(t)&\geq y^{\ep}(t_1)+k_0\gamma^{\frac{2}{n}}\int\limits_{t_1}^{t}\!\int\limits_0^{\infty}\varphi\chie W^{\ep}\intd s\intd t-\frac{1}{2}\int\limits_{t_1}^{t}\!\int\limits_0^{\infty}\varphi_s(\chie)^2(W^{\ep})^2\intd s\intd t\nonumber\\
&+n\int\limits_{t_1}^{t}\!\int\limits_0^{\infty} F\varphi_s(1-\chie)\chie W^{\ep}\intd s\intd t+n\int\limits_{t_1}^{t}\!\int\limits_0^{\infty} F_s\varphi(1-\chie)\chie W^{\ep}\intd s\intd t\nonumber\\
&-c_2\varepsilon^{\beta-\delta}(t-t_1)\ \text{ for all }t\in(t_1,t_0+\eta)\text{ and }0<\varepsilon<\min\left\lbrace1,\frac{\xi}{\gamma}\right\rbrace.
\end{align}
In order to take $\varepsilon\searrow0$ we recall the definition of $\varphi$ in \eqref{phidef} to see that
\begin{align*}
|\varphi_s(s)|\leq\left.\begin{cases}
\frac{a\delta}{\gamma^\delta}s^{-(\delta+1)}\leq \frac{a\delta}{\gamma^\delta}e^\xi e^{-\gamma s}\left(1+s^{-(\delta+1)}\right)&\text{ if }s< \frac{\xi}{\gamma}\\
\gamma e^{-\gamma s}\leq\gamma\left(1+s^{-(\delta+1)}\right)e^{-\gamma s}&\text{ if }s> \frac{\xi}{\gamma}
\end{cases}\right\rbrace\leq c_3\left(1+s^{-(\delta+1)}\right)e^{-\gamma s},
\end{align*}
for some $c_3>0$ and $s>0$. And similarly
\begin{align*}
|\varphi(s)|\leq c_4\left(1+s^{-\delta}\right)e^{-\gamma s},
\end{align*}
for some $c_4>0$ in $(0,\infty)$. Combining these inequalities with our assumption $\eqref{Wscontra}$ and the definitions of $F$ and $F_s$ in \eqref{Fdef} and \eqref{Fsdef}, respectively, we observe that
\begin{align*}
\left\vert F(s)\varphi_s(s)(1-\chie\!(s))\!\chie\!(s) W^{\ep}\!(s,t)\right\vert\leq\left\vert F(R+\rho)\varphi_s(s)W(s,t)\right\vert\leq c_5(s^\beta\!+s^{\beta-\delta-1})e^{-\gamma s}
\end{align*}
for some $c_5>0$ a.e in $(0,\infty)$, as well as
\begin{align*}
\left\vert F_s(s)\varphi(s)(1-\chie\!(s))\chie\!(s)W^{\ep}\!(s,t)\right\vert\leq\left\vert F_s(s)\varphi(s)W(s,t)\right\vert\leq c_6(s^{\beta-\frac{\alpha}{n}}+s^{\beta-\frac{\alpha}{n}-\delta})e^{-\gamma s}
\end{align*}
for some $c_6>0$ and all $s>0$, holds independent of $\varepsilon$. Because of $\beta>\delta$ and $n>\alpha$ the integrals $\smallint\limits_{t_1}^{t}\!\smallint\limits_0^{\infty} (s^\beta+s^{\beta-\delta-1})e^{-\gamma s}\intd s\intd t$ and $\smallint\limits_{t_1}^{t}\!\smallint\limits_0^{\infty}(1+s^{\beta-\delta-\frac{\alpha}{n}})e^{-\gamma s}\intd s\intd t$ converge, so that an additional application of the dominated convergence theorem shows
\begin{align*}
\lim_{\varepsilon\searrow0}\int\limits_{t_1}^{t}\!\int\limits_0^{\infty}F\varphi_s(1-\chie)\chie W^{\ep}\intd s\intd t= \int\limits_{t_1}^{t}\int\limits_0^{\infty}\lim_{\varepsilon\searrow0} F\varphi_s(1-\chie)\chie W^{\ep}\intd s\intd t=0
\end{align*}
and
\begin{align*}
\lim_{\varepsilon\searrow0}\int\limits_{t_1}^{t}\!\int\limits_0^{\infty}F_s\varphi(1-\chie)\chie W^{\ep}\intd s\intd t=0
\end{align*}
for all $t\in(t_1,t_0+\eta)$. For the two remaining integral terms in \eqref{blowupproofeq10} we make use of the monotonicity of $\chie$ and $W^{\ep}$ with respect to $\varepsilon$ and apply the monotone convergence theorem to conclude that $y(t):=\smallint\limits_0^\infty\varphi W\intd s$ satisfies
\begin{align}\label{blowupproofeq11}
\owncount
y(t)&\geq y(t_1)+k_0\gamma^{\frac{2}{n}}\int\limits_{t_1}^{t}\!\int\limits_0^{\infty}\varphi W\intd s\intd t-\frac{1}{2}\int\limits_{t_1}^{t}\!\int\limits_0^{\infty}\varphi_s W^2\intd s\intd t\nonumber\\
&= y(t_1)+k_0\gamma^{\frac{2}{n}}\int\limits_{t_1}^{t}\!\int\limits_0^{\infty}\varphi W\intd s\intd t+\frac{1}{2}\int\limits_{t_1}^{t}\!\int\limits_0^{\infty}|\varphi_s| W^2\intd s\intd t\ \text{ for all }t\in(t_1,t_0+\eta).
\end{align}
Using H\"older's inequality and the definition of $y(t)$ we see that
\begin{align*}
y^2(t)&=\left(\int\limits_0^\infty\varphi W\intd s\right)^2\!
\leq
\left(\int\limits_0^\infty\frac{\varphi^2}{\vert\varphi_s\vert}\intd s\right)\left(\int\limits_0^\infty\vert\varphi_s\vert W^2\intd s\right)\ \text{ for all }t>0,
\end{align*}
which in turn by Lemma \ref{lemtestfunc} \eqref{phiint} implies
\begin{align*}
y^2(t)\leq\frac{K_0}{\gamma^2}\int\limits_0^\infty\vert\varphi_s\vert W^2\intd s\ \text{ for all }t>0.
\end{align*}
Combination with \eqref{blowupproofeq11} therefore yields
\begin{align*}
y(t)\geq y(t_1)+k_0\gamma^{\frac{2}{n}}\int\limits_{t_1}^{t}y(\tau)\intd\tau+\frac{\gamma^2}{2K_0}\int\limits_{t_1}^t y^2(\tau)\intd\tau\ \text{ for all }t\in(t_1,t_0+\eta).
\end{align*}
In view of Lemma \ref{grönvar} $y(t)$ thus satisfies the inequality $y(t)\geq z(t)$ for all $t\in(t_1,t_0+\eta)$, where $z$ is the solution of
\begin{align*}
\left\{\begin{array}{c@{\,}l@{\,}c}
z'&=Az+Bz^2,\ &t>t_1,\\
z(t_1)&=y(t_1)&
        \end{array}\right.
\end{align*}
with $A=k_0\gamma^{\frac{2}{n}}$ and $B=\frac{\gamma^2}{2K_0}$. For $y(t_1)>0$ this Bernoulli-type initial-value problem has the explicit solution
\begin{align*}
z(t)=\frac{1}{\left(\frac{1}{y(t_1)}+\frac{B}{A}\right)e^{-A(t-t_1)}-\frac{B}{A}},\ t\in(t_1,t_1+T),
\end{align*}
with maximal existence time determined by $T=\frac{1}{A}\log\left(1+\frac{A}{By(t_1)}\right)$. Next, we utilize $\kappa\gamma^{\frac{2-n}{n}}>\frac{\xi}{\gamma}$ from \eqref{gammaineqR} and the fact $W_s\geq0$, to estimate $y(t_1)$ from below with
\begin{align*}
y(t_1)&=\int\limits_0^\infty\varphi(s)W(s,t_1)\intd s\geq c_{\gamma}\hspace*{-5pt}\int\limits_{\kappa\gamma^{\frac{2-n}{n}}}^\infty\hspace*{-5pt}\varphi(s)\intd s\geq c_{\gamma}\hspace*{-7pt}\int\limits_{\kappa\gamma^{\frac{2-n}{n}}}^\infty\hspace*{-4pt} e^{-\gamma s}\intd s=\frac{c_{\gamma}}{\gamma}e^{-\kappa\gamma^{\frac{2}{n}}},
\end{align*}
where we set $c_{\gamma}:=W(\kappa\gamma^{\frac{2-n}{n}},t_1)$. Accordingly, recalling inequality \eqref{gammaineqe} from Lemma \ref{lempara} we have
\begin{align*}
T&\leq\frac{1}{k_0\gamma^{\frac{2}{n}}}\log\left(1+\frac{k_0\gamma^{\frac{2}{n}}}{\frac{\gamma^2c_{\gamma}}{2K_0\gamma}e^{-\kappa\gamma^{\frac{2}{n}}}}\right)=\frac{1}{k_0\gamma^{\frac{2}{n}}}\log\left(1+\frac{2k_0 K_0}{c_{\gamma}\gamma^{\frac{n-2}{n}}}e^{\kappa\gamma^{\frac{2}{n}}}\right)\\
&\leq\frac{1}{k_0\gamma^{\frac{2}{n}}}\log\left(e^{2\kappa\gamma^{\frac{2}{n}}}\right)=\frac{2\kappa}{k_0}
\end{align*}
by definitions of $A$ and $B$. But this means, see \eqref{kappaineq}, that $T\leq\frac{\eta}{4}$. Thus $y$ blows up before or at $t=t_1+\frac{\eta}{4}=t_0+\frac{3\eta}{4}$, which -- since $\varphi$ is integrable because of $\delta<1$ -- is a contradiction to
\begin{align*}
y(t)=\int\limits_0^\infty\varphi(s) W(s,t)\intd s\leq\frac{\mu n}{\vert S_{n-1}\vert}\int\limits_0^\infty\varphi(s)\intd s\leq c_6\ \forall t>0.
\end{align*}
Hence our assumption in \eqref{Wscontra} must have been false, which completes the proof of \eqref{wsblow}. To verify \eqref{wsblow2} we take $\beta=1$ and conclude that $W(s,t)$ cannot be Lipschitz continuous on $(0,\infty)\times(t_0,t_0+\eta)$ for each $t_0\geq0$ and every $\eta>0$, which then implies the asserted unboundedness of $W_s$.
\end{bew}
Having this result for the derivative of $W$ at hand, we can now show our main result, which corresponds in part to \cite[Theorem 0.1]{TW11}.
\begin{proof}[\textbf{Proof of Theorem \ref{theoblowup}:}]
The assumption $u_0\equiv c_0$ on $\overline{B_1(0)}$ for some $c_0>0$, implies $W_0\in C^0([0,1])\cap C^2((0,1))$, $W_0(s)=c_0 s$ on $[0,1]$, $W_{0s}=c_0>0$ on $(0,1)$ and $W_{0ss}=0$ on $(0,1)$. This allows us to choose one of our generalized test functions $\varphi$ fulfilling the important inequalities \eqref{phiodi} and \eqref{phiint} of Lemma \ref{lemtestfunc}, with parameters chosen as in Lemma \ref{lempara}.
Using this test function as show in by Proposition \ref{wspropo}, we obtain $\|W_s\|_{\Lo^\infty((0,\infty)\times(t_0,t_0+\eta))}=\infty$ for all $\eta>0$ (see \eqref{wsblow2}). Now we recall the measure-valued reconstruction identity $u(x,t)=W_s(\vert x\vert^n\!,t)+\frac{\vert S_{n-1}\vert}{n} W(0+,t)\delta(x)$ (shown in \eqref{backtrafo}) to verify that the radial weak solution of \eqref{KS0f}, in the sense of Definition \ref{weakrad}, blows up immediately at $x=0$.
\end{proof}


\end{document}